\documentclass[11pt]{article}

\usepackage[T2A]{fontenc}
\usepackage[cp1251]{inputenc}
\usepackage{amsthm}
\usepackage{amsfonts}
\usepackage{eufrak}
\usepackage{amssymb}
\usepackage{amsmath}
\usepackage{cite}

\DeclareMathOperator{\Ker}{Ker}
\begin{document}
\newtheoremstyle{mytheorem}
  {\topsep}   
  {\topsep}   
  {\itshape}  
  {}       
  {\bfseries} 
  { }         
  {5pt plus 0pt minus 1pt} 
  { }          
\newtheoremstyle{myremark}
  {\topsep}   
  {\topsep}   
  {\upshape}  
  {}       
  {\bfseries} 
  {  }         
  {5pt plus 1pt minus 1pt} 
  { }          
\theoremstyle{mytheorem}
\newtheorem{theorem}{Theorem}[section]
 \newtheorem{theorema}{Theorem}
 \newtheorem*{heyde1*}{Theorem A}
 \newtheorem*{heyde2*}{Theorem B}
 \newtheorem*{heyde3*}{Теорема C}
 \newtheorem{lemma}[theorem]{Lemma}
\newtheorem{corollary}[theorem]{Corollary}
\theoremstyle{myremark}
\newtheorem{remark}[theorem]{Remark}
\noindent This article is accepted for publishing in

\noindent Journal of Mathematical Analysis and Applications

\vskip 1 cm

\centerline{\bf On a characterization theorem in the space $\mathbb{R}^n$}

\bigskip

\centerline{\textbf{G.M. Feldman}}

\bigskip

\makebox[20mm]{ }\parbox{125mm}{ \small By Heyde's theorem,  the class of Gaussian distributions on the real line is characterized by the symmetry of the conditional distribution of one linear form of independent random variables given another.   We prove an  analogue of this theorem for two independent random vectors taking values in the space $\mathbb{R}^n$. The obtained class of distributions consists of    convolutions of Gaussian distributions and  a distribution supported in a subspace, which is determined by   coefficients of the linear forms.}

\bigskip


\bigskip


\section{ Introduction}

By the well-known Skitovich--Darmois theorem,  Gaussian distributions on the real are characterized by the independence of two linear forms with nonzero coefficients of independent random variables. This result was generalized by S.G.~Ghurye  and I.~Olkin to independent random vectors taking values in the space $\mathbb{R}^n$. Coefficients of the linear forms in this case are invertible linear operators in $\mathbb{R}^n$ (\!\!\cite{GhO}, see also \cite[\S 3.2]{KaLiRa}). A theorem similar to the Skitovich--Darmois theorem was proved by C.C.~Heyde.
By   Heyde's  theorem, Gaussian distributions on the real line are characterized by the   symmetry of the conditional distribution of one linear form   given another. Coefficients of the linear forms are nonzero real numbers (\!\!\cite{He}, see also \cite[\S 13.4]{KaLiRa}). For two independent random variables 
Heyde's theorem states the following. Let $\xi_1$ and $\xi_2$ be
independent random variables with distributions
$\mu_1$ and $\mu_2$. Let $a_j$ and $b_j$ be nonzero numbers such that
$b_1a_1^{-1}+b_2a_2^{-1}\ne 0$. If the conditional  distribution 
of linear form $L_2=b_1\xi_1+b_2\xi_2$ given $L_1=a_1\xi_1+a_2\xi_2$
 is symmetric, then  $\mu_j$  are Gaussian distributions.
 It is easy to see that in studying the possible distributions  
$\mu_j$  we can assume without loss of generality that $L_1=\xi_1+\xi_2$, $L_2=\xi_1+a\xi_2$, where $a\ne 0$. Hence the following statement holds.

\begin{heyde1*}  Let $\xi_1$ and $\xi_2$ be
independent random variables with distributions
$\mu_1$ and $\mu_2$. Let $a\ne 0$ and $a\ne -1$. 
  If the conditional  distribution of the linear form 
  $L_2 = \xi_1 +a\xi_2$ given $L_1 = \xi_1 + \xi_2$ is symmetric,
then  $\mu_j$  are Gaussian distributions.
\end{heyde1*}
 
A number of works have been devoted to 
analogues of   Heyde's theorem for various locally compact Abelian groups (see e.g. \cite{Fe2, MiFe1, Fe4, Fe20bb, My2, Fe6, FeTVP1,  M2013, M2020, JFAA2021, POTA2022, Fe2015a}). In so doing  coefficients of the linear forms are  topological automorphisms of the group. At the same time  
 Heyde's theorem for  the space $\mathbb{R}^n$ was not specially studied. The space $\mathbb{R}^n$ was considered only as a special case of  a locally compact Abelian group. In particular, the following theorem was proved in  \cite{FeTVP1} for arbitrary  locally compact Abelian groups containing no nonzero elements of finite order. For the space $\mathbb{R}^n$ this theorem can be formulated as follows.
\begin{heyde2*}\label{A}   Let    $\alpha$ be an invertible linear operator in the space $\mathbb{R}^n$ such that $I+\alpha$ is also an invertible  operator.  Let $\xi_1$ and $\xi_2$ be
independent random vectors with values in $\mathbb{R}^n$ and distributions
$\mu_1$ and $\mu_2$. If the   conditional  distribution of the linear form  $L_1 = \xi_1 + \xi_2$ given $L_2 = \xi_1 +
\alpha\xi_2$  is symmetric, then  $\mu_j$  are Gaussian distributions.
\end{heyde2*}
Obviously, Theorem A is a particular case of Theorem B for   $n=1$.
The purpose of this note is to prove  Heyde's theorem for the space     $\mathbb{R}^n$ without any restrictions on an   invertible linear 
operator $\alpha$. In other words, we want to get   a full  description of the possible distributions  $\mu_j$ in Theorem B    in the case when  the  
operator $I+\alpha$ need not be invertible.  It turns out that if the   
operator $I+\alpha$ is not   invertible, then $\mu_j$ are convolutions of Gaussian  distributions and a distribution supported in a subspace which is determined by the  operator $\alpha$.

Denote by $x=(x_1, x_2, \dots, x_n)$, $x_j\in \mathbb{R}$, elements of the space  $\mathbb{R}^n$.    If $x, y\in \mathbb{R}^n$, then put   $$\langle x, y\rangle=\sum_{j=1}^nx_jy_j.$$
Denote by $\mathbb{C}$ the complex plane. We will also use the notation $\langle x, y\rangle$ in the case when   $x, y\in \mathbb{C}^n$. Let $H$ be a subspace of   $\mathbb{R}^n$. Denote by $$A(\mathbb{R}^n, H)=\{x\in \mathbb{R}^n: \langle x, y\rangle=0 \ \mbox{for all} \ y\in H\}$$ the annihilator of $H$. 

Let $\alpha$ be a linear operator in    $\mathbb{R}^n$. Denote by $\widetilde\alpha$   its adjoint  operator. If   a subspace $G$ of   $\mathbb{R}^n$  is invariant with respect to   $\alpha$, then denote by $\alpha{_G}$ the restriction of  $\alpha$ to   $G$. Denote by $I$ the identity operator. Denote by $\|x\|$ a norm of a vector $x\in \mathbb{R}^n$  
and by $\|\alpha\|$ the  norm of the operator $\alpha$.

Let $P(y)$ be an arbitrary function on $\mathbb{R}^n$  and let $h\in \mathbb{R}^n$.
 Denote by $\Delta_h$ the finite difference operator
$$\Delta_h P(y)=P(y+h)-P(y), \quad y\in \mathbb{R}^n.$$
Recall that a continuous function    $P(y)$ is a polynomial in some neighbourhood of zero in   $\mathbb{R}^n$ if and only if for a nonnegative integer    $m$ the function $P(y)$ satisfies the equation
\begin{equation}\label{18.05.1}
\Delta_{h}^{m+1}P(y)=0,
\end{equation}
for all $y$ and $h$ in a neighborhood of zero in  $\mathbb{R}^n$. Moreover, the minimum $m$  at which
 (\ref{18.05.1}) is satisfied coincides with the degree of the poynomial $P(y)$.

Let  $\mu$ be a probability distribution on   $\mathbb{R}^n$. Denote by $\hat\mu(y)$, $y\in \mathbb{R}^n$, the characteristic function of    $\mu$.    Define the distribution   $\bar\mu$  by the formula
$\bar\mu(B) =\mu(-B)$ for each Borel subset $B$ in $\mathbb{R}^n$. Then $\hat{\bar{\mu}}(y)=\overline{\hat\mu(y)}$. Denote by $E_x$
the degenerate distribution
 concentrated at a vector $x\in \mathbb{R}^n$.

\section{ Main theorem}

The main result of this paper is the proof of the following theorem.
\begin{theorem}\label{th1}   Let    $\alpha$ be an invertible linear operator in the space $\mathbb{R}^n$.  Put $K={\Ker}~(I+\alpha)$.  Let $\xi_1$ and $\xi_2$ be
independent random vectors with values in $\mathbb{R}^n$ and distributions
$\mu_1$ and $\mu_2$. Assume that   the conditional  distribution of the linear form  $L_2 = \xi_1 + \alpha\xi_2$ given $L_1 = \xi_1 + \xi_2$  is symmetric.
Then there exists an  $\alpha$-invariant subspace $G$  such that $\mu_j$ are shifts of convolutions of  symmetric Gaussian distributions supported in   $G$   and a distribution supported in  $K$. Moreover, $K\cap G=\{0\}$.
\end{theorem}
It is clear that if in Theorem \ref{th1} $I+\alpha$ is an invertible operator, i.e. $K=\{0\}$, then $\mu_j$ are   Gaussian distributions. Thus, Theorem B
follows from Theorem \ref{th1}.

To prove Theorem \ref{th1} we need a series of lemmas.
The following lemma holds for independent random variables taking values 
in an arbitrary locally compact Abelian group. We formulate it for 
the space $\mathbb{R}^n$.
\begin{lemma}\label{lem1}{\rm(\!\!\cite[Lemma 16.1]{Fe5})}  Let    $\alpha$ be an invertible linear operator in the space $\mathbb{R}^n$. Let $\xi_1$ and $\xi_2$ be
independent random vectors with values in $\mathbb{R}^n$ and distributions
$\mu_1$ and $\mu_2$. The conditional  distribution of the linear form  $L_2 = \xi_1 + \alpha\xi_2$ given $L_1 = \xi_1 + \xi_2$  is symmetric if and only if the characteristic functions   $\hat\mu_j(y)$  satisfy the equation
\begin{equation}\label{11.04.1}
\hat\mu_1(u+v )\hat\mu_2(u+\widetilde\alpha v )=
\hat\mu_1(u-v )\hat\mu_2(u-\widetilde\alpha v), \quad u, v \in \mathbb{R}^n.
\end{equation}
\end{lemma}
Equation (\ref{11.04.1}) is called the Heyde  functional 
 equation. Due to Lemma \ref{lem1}, the proof of  Theorem \ref{th1} reduces to
the description of solutions of equation (\ref{11.04.1}) in the class of continuous normalized positive definite functions. Note that the proof of the Ghurye--Olkin  theorem, mentioned in the introduction, for two independent random vectors reduces to solving the Skitovich--Darmois functional equation
\begin{equation}\label{new04.12.22}
\hat\mu_1(u+v )\hat\mu_2(u+\widetilde\alpha v )=
\hat\mu_1(u)\hat\mu_2(u)\hat\mu_1(v)\hat\mu_2(\widetilde\alpha v), \quad u, v \in \mathbb{R}^n,
\end{equation}
in the class of continuous normalized positive definite functions. Unlike  
 equation (\ref{11.04.1}), all  solutions of equation (\ref{new04.12.22}) 
 in the class of continuous normalized positive definite functions 
 are characteristic functions of Gaussian distributions in the space $\mathbb{R}^n$.

Note that based on characterization  of polynomials as the solutions set of
some functional equations, J. M. Almira in \cite{Al1} proposed a new approach
to solving the Skitovich--Darmois functional equation 
in the space $\mathbb{R}^n$ for $m\ge 2$ 
 functions. Then, using the fact that Aichinger’s equation characterizes polynomial functions (see \cite{Ai} by E. Aichinger and J. Moosbauer), 
 J. M. Almira in \cite{Al2} studied the solutions of the Skitovich--Darmois functional equation 
 on an arbitrary Abelian group.

\medskip   

We will need the following easily verified statement, which we formulate as a lemma (see e.g.  \cite[Lemma 6.9]{Fe9}).
\begin{lemma}\label{lem5}
Let $\mathbb{R}^n=\mathbb{R}^p\times \mathbb{R}^q$  and let
$\mu$ be a distribution on the space   $\mathbb{R}^n$ with the characteristic function   $\hat\mu(s_1, s_2)$, $s_1\in \mathbb{R}^p$, $s_2\in \mathbb{R}^q$. Assume that the function $\hat\mu(0, s_2)$, $s_2\in \mathbb{R}^q$, is extended to   $\mathbb{C}^q$ as an entire function in
$s_2$. Put $B_r=\{s_2=(s_{21}, s_{22},\dots, s_{2q})\in \mathbb{C}^q: |s_{2j}|\le r, \ j=1, 2, \dots, q\}$. Then for each fixed    $s_1\in \mathbb{R}^p$ the function    $\hat\mu(s_1, s_2)$, $s_2\in \mathbb{R}^q$, is also extended to   $\mathbb{C}^q$
as an entire function in  $s_2$,    and  for each $s_1\in \mathbb{R}^p$ the inequality
\begin{equation}\label{16.05.6}
\max_{ s_2\in B_r}|\hat\mu(s_1, s_2)|\le \max_{s_2\in B_r}|\hat\mu(0, s_2)|
\end{equation}
holds.
\end{lemma}
The following lemma holds for independent random variables taking values in an arbitrary locally compact Abelian group containing no elements of order 2.    It follows directly from Lemma  \ref{lem1}. We formulate it for the space  $\mathbb{R}^n$.
\begin{lemma}\label{lem3}
Let $\xi_1$ and $\xi_2$ be
independent random vectors with values in the space $\mathbb{R}^n$ and distributions
$\mu_1$ and $\mu_2$. The conditional  distribution of the linear form  $L_2 = \xi_1 -\xi_2$ given $L_1 = \xi_1 + \xi_2$  is symmetric if and only if  $\mu_1=\mu_2$.
\end{lemma}
The following lemma is crucial in the proof of Theorem \ref{th1}. It describes the possible distributions  $\mu_j$ in Theorem B in the case when in a suitable basis  a Jordan cell with the  eigenvalue $\lambda=-1$ corresponds to a linear operator $\alpha$, i.e.
$$\alpha=\alpha_n=\left(
\begin{array}{ccccc}
-1 & 1 & 0  & \ldots & 0\\
0  & -1 & 1 & \ldots & 0\\
\ldots & \ldots & \ldots &\ldots   &\ldots\\
0 & 0 & 0 & \ldots & 1 \\
0 & 0 & 0 & \ldots & -1
\end{array}
\right).
$$
\begin{lemma}\label{lem13}   Let    $\alpha_n$  be the invertible linear operator in the space $\mathbb{R}^n$, $n\ge 2$,   of the form
\begin{equation}\label{16.05.1}
\alpha_n(x_1, x_2, \dots, x_n)=(-x_1+x_2, -x_2+x_3, \dots,-x_{n-1}+ x_n, -x_n), \quad x_j\in \mathbb{R}.
\end{equation}
Let $\xi_1$ and $\xi_2$ be
independent random vectors with values in $\mathbb{R}^n$ and distributions
$\mu_1$ and $\mu_2$. Assume that the conditional  distribution of the 
linear form  $L_2 = \xi_1 + \alpha_n\xi_2$ given $L_1 = \xi_1 + \xi_2$  is symmetric.   
Put $K={\Ker}~(I+\alpha_n)$. Then we can replace the distributions $\mu_j$ by their shifts $\tau_j$ in such a way that  $\tau_1 =\tau_2$, the distribution $\tau_j$ is supported in the subspace $K$, and if $\eta_j$ are independent identically distributed random vectors with values in the space $\mathbb{R}^n$ and distribution  $\tau_j$, then the conditional  distribution of the linear form
$M_2 = \eta_1 +\alpha_n\eta_2$ given $M_1 = \eta_1 + \eta_2$  is symmetric.  \end{lemma}
{\bf Proof.} Note that $K=\{(x_1, 0, \dots, 0)\in \mathbb{R}^n:x_1\in \mathbb{R}\}$. We divide the proof of the lemma into three steps.

1. By Lemma \ref{lem1}, the characteristic functions  $\hat\mu_j(y)$   satisfy equation  (\ref{11.04.1}). Put $\nu_j=\mu_j*\bar\mu_j$. Then $\hat\nu_j(y)=|\hat\mu_j(y)|^2\ge 0$. Since $\hat\nu_j(0)=1$, we can choose $\varepsilon >0$ in such a way that the inequalities $\hat\nu_j(y)>0$, $j=1, 2$, are valid for  all $\|y\|<\varepsilon$.   Obviously, the characteristic functions $\hat\nu_j(y)$ also  satisfy equation  (\ref{11.04.1}). Put
$P_j(y)=\log\hat \nu_j(y)$, where $\|y\|<\varepsilon$, $j=1, 2$. Let $\|\alpha_n\|=M$. It is obvious that $M\ge 1$. Put $\delta=\varepsilon/8M$  and
$$H=(I+\widetilde\alpha_n)(\mathbb{R}^n).$$
We prove that $P_j(y)$ are polynomials of degree at most 2 in some neighbourhood of zero in the subspace $H$.

 It follows from equation   (\ref{11.04.1}) for the characteristic functions  $\hat\nu_j(y)$ that the functions $P_j(y)$ satisfy the equation
\begin{equation}\label{14.05.1}
P_1(u+v)+P_2(u+\widetilde\alpha_n v)-P_1(u-v)-P_2(u-\widetilde\alpha_n v)=0, \quad u, v \in \mathbb{R}^n, \ \|u\|<\delta, \ \|v\|<\delta.
\end{equation}
Equation (\ref{14.05.1})  arises in the study of    Heyde's theorem on various locally compact Abelian groups. To solve equation (\ref{14.05.1}) we use the finite difference method. This is a standard reasoning
(see e.g. \cite{Fe6}, \cite{FeTVP1}).

Take an arbitrary vector $k_1$ in the space   $\mathbb{R}^n$ such that $\|k_1\|<\delta$. Put $h_1
=\widetilde\alpha_n k_1$.   Replacing in (\ref{14.05.1})   $u$ by $u+h_1$  and $v$ by $v+k_1$ and  subtracting   equation
(\ref{14.05.1}) from the obtained equation, we get
\begin{equation}\label{14.05.2}
\Delta_{l_{11}}P_1(u+v) +
\Delta_{l_{12}}P_2(u+\widetilde\alpha_n v) -
\Delta_{l_{13}}P_1(u-v)=0, \quad u, v \in \mathbb{R}^n, \ \|u\|<\delta, \ \|v\|<\delta,
\end{equation}
where $l_{11}= (I+\widetilde\alpha_n)k_1$, $l_{12}=2 \widetilde\alpha_n k_1$, $l_{13}= (\widetilde\alpha_n-I)k_1.$
Take an arbitrary vector $k_2$ in  the space   $\mathbb{R}^n$ such that $\|k_2\|<\delta$. Put $h_2
=k_2$.   Replacing in (\ref{14.05.2})  $u$ by $u+h_2$  and $v$ by $v+k_2$ and  subtracting   equation
(\ref{14.05.2}) from the obtained equation, we find
\begin{equation}\label{14.05.3}
\Delta_{l_{21}}\Delta_{l_{11}}P_1(u+v) +
\Delta_{l_{22}}\Delta_{l_{12}}P_2(u+\widetilde\alpha_n v) =
0, \quad u, v \in \mathbb{R}^n, \ \|u\|<\delta, \ \|v\|<\delta,
\end{equation}
where $l_{21}=2k_2$, $l_{22}=
(I+\widetilde\alpha_n)k_2.$
Take an arbitrary vector $k_3$ in the space   $\mathbb{R}^n$ such that $\|k_3\|<\delta$. Put $h_3
=-\widetilde\alpha_n k_3$.   Replacing in (\ref{14.05.3})   $u$ by $u+h_3$  and $v$ by $v+k_3$ and  subtracting   equation
(\ref{14.05.3}) from the obtained equation, we get
\begin{equation}\label{14.05.4}
\Delta_{l_{31}}\Delta_{l_{21}}\Delta_{l_{11}}P_1(u+v)
= 0, \ \ u, v \in \mathbb{R}^n, \ \|u\|<\delta, \ \|v\|<\delta,
\end{equation}
where $l_{31}= (I-\widetilde\alpha_n)k_3.$ Substituting
    $v=0$ into (\ref{14.05.4})  we find that
\begin{equation}\label{14.05.5}
\Delta_{l_{31}}\Delta_{l_{21}}\Delta_{l_{11}}P_1(u) = 0, \quad u \in \mathbb{R}^n, \ \|u\|<\delta.
\end{equation}
Since multiplication by   2 and $I-\widetilde\alpha_n$ are invertible linear operators in the space   $\mathbb{R}^n$  and $k_j$ are arbitrary vectors in   $\mathbb{R}^n$ such that $\|k_j\|<\delta$, it follows from (\ref{14.05.5})
   that the function $P_1(y)$ satisfies the equation
\begin{equation}\label{14.05.6}
\Delta_h^3P_1(y)=0, \ \ y, h \in H,  \  \|y\|<\delta, \ \|h\|<\delta.
\end{equation}
Arguing similarly, we exclude the function $P_1(y)$ from equation   (\ref{14.05.3}) and obtain that the function    $P_2(y)$ also satisfies  equation    (\ref{14.05.6}). It follows from this that $P_j(y)$ are polynomials of degree at most 2 in some neighbourhood of zero in the subspace $H$.

2. Assume $n=2$ and prove that $\hat\nu_1(y)=\hat\nu_2(y)=1$ for all $y\in H$. Since $\alpha_2(x_1, x_2)=(-x_1+x_2, -x_2)$, $x_j\in \mathbb{R}$, we have $\widetilde\alpha_2(y_1, y_2)=(-y_1, y_1-y_2)$, $y_j\in \mathbb{R}$, and then
$H=\{(0, y_2)\in \mathbb{R}^2: y_2\in \mathbb{R}\}$.  Since $P_j(y)$ are polynomials of the degree at most 2 in some neighbourhood of zero in the subspace $H$   and $P_j(y)=\log\hat \nu_j(y)$ for all $\|y\|<\delta$,  $y\in H$,  this implies that the   characteristic functions   $\hat \nu_j(0, y_2)$ coincide with the characteristic functions of some symmetric Gaussian distributions in some neighbourhood of zero on the real line. It is well known that then $\hat \nu_j(0, y_2)$, $y_2\in \mathbb{R}$, are the characteristic functions of some symmetric Gaussian distributions (\!\!\cite[\S 1.2]{KaLiRa}).  Thus, we have
\begin{equation}\label{14.05.7}
P_j(0, y_2)=-\sigma_j y_2^2, \quad  y_2\in \mathbb{R},
\end{equation}
where $\sigma_j\ge 0$, $j=1, 2$.

Since $\widetilde\alpha_2(y_1, y_2)=(-y_1, y_1-y_2)$, $y_j\in \mathbb{R}$, equation (\ref{11.04.1}) for the functions $\hat\nu_j(y)$ takes the form
\begin{multline}\label{14.05.8}
\hat\nu_1(u_1+v_1, u_2+v_2)\hat\nu_2(u_1- v_1, u_2+v_1-v_2)\\=
\hat\nu_1(u_1-v_1, u_2-v_2)\hat\nu_2(u_1+ v_1, u_2-v_1+v_2), \quad  u_j, v_j\in \mathbb{R}.
\end{multline}
By Lemma \ref{lem5}, it follows from (\ref{14.05.7}) that for each fixed $y_1\in \mathbb{R}$ the functions $\hat\nu_j(y_1, y_2)$ can be extended to the   complex plane  $\mathbb{C}$   as entire functions in $y_2$.  It is obvious that equation  (\ref{14.05.8}) remains valid  for all $u_1, v_1  \in \mathbb{R}$ and $u_2, v_2  \in \mathbb{C}$.

We will verify that if   $\hat\nu_j(y_1, 0)>0$ for  all $|y_1|<\delta$, $j=1, 2$, then the functions $\hat\nu_j(y_1, y_2)$ do  not vanish for all    $|y_1|<\delta$ and      $y_2  \in \mathbb{C}$. Assume the contrary, let $\hat\nu_1(\widetilde y_1, \widetilde y_2)=0$, where $|\widetilde y_1|<\delta$, $\widetilde y_2 \in \mathbb{C}$. Take arbitrary numbers  $u_1   \in \mathbb{R}$ and    $u_2  \in \mathbb{C}$. Put   $v_1=\widetilde y_1-u_1$, $v_2=\widetilde y_2-u_2$ and substitute  $u_j$ and $v_j$ into equation  (\ref{14.05.8}). Then the left-hand side of equation (\ref{14.05.8}) is equal to zero. Hence  the equality
\begin{equation}\label{14.05.10}
\hat\nu_1(2u_1-\widetilde y_1, 2u_2-\widetilde y_2)\hat\nu_2(\widetilde y_1, \widetilde y_2-\widetilde y_1+u_1)=0
\end{equation}
is satisfied  for all $u_j$. Since $\hat\nu_2(\widetilde y_1, y_2)$ is not identically zero an entire function in $y_2$, take $u_1=\widetilde u_1$ in such a way that $\hat\nu_2(\widetilde y_1, \widetilde y_2-\widetilde y_1+\widetilde u_1)\ne 0$ and $|2\widetilde u_1-\widetilde y_1|<\delta$. Obviously, we can do it. Then it follows from (\ref{14.05.10}) that $\hat\nu_1(2\widetilde u_1-\widetilde y_1, 2u_2-\widetilde y_2)=0$ for  all $u_2\in \mathbb{C}$, that  contradicts the fact that $\hat\nu_1(2\widetilde u_1-\widetilde y_1, y_2)$ is a not identically zero entire function in  $y_2$. Arguing similarly it is easy to make sure that the function  $\hat\nu_2(y_1, y_2)$ does not vanish for all $|y_1|<\delta$, $y_2  \in \mathbb{C}$.

Taking into account (\ref{14.05.7}), by Lemma \ref{lem5}, it follows from inequality   (\ref{16.05.6}) that           $\hat\nu_j(y_1, y_2)$ are entire functions in   $y_2$  of order at most 2 for all $|y_1|<\delta$. Since the functions $\hat\nu_j(y_1, y_2)$ for all $|y_1|<\delta$, $y_2  \in \mathbb{C}$ do not vanish,    Hadamard's theorem on the representation of entire functions of finite order implies that the representation
\begin{equation}\label{15.05.1}
\hat\nu_j(y_1, y_2)=\exp\{a_j(y_1) y_2^2+b_j(y_1)y_2+c_j(y_1)\}, \quad |y_1|<\delta, \  y_2  \in \mathbb{C},
\end{equation}
is valid, where  $a_j(y_1), b_j(y_1), c_j(y_1)\in \mathbb{C}$, $j=1, 2$.
Substituting $u_1=v_1=0$, $u_2=v_2=y_2/2$ into  (\ref{14.05.8}), we get
\begin{equation}\label{14.05.9}
\hat\nu_1(0, y_2)=\hat\nu_2(0, y_2), \ \  y_2\in \mathbb{R}.
\end{equation}
It follows from (\ref{14.05.7}) and (\ref{14.05.9})  that $\sigma_1=\sigma_2=\sigma$.  Thus, (\ref{14.05.7}) and (\ref{15.05.1}) imply that $a_1(0)=a_2(0)=-\sigma$, $b_1(0)=b_2(0)=0$. Obviously, without loss of generality, we can assume that    $c_1(0)=c_2(0)=0$.

Substituting $u_1=v_1=y_1/2$ into (\ref{14.05.8})   and taking into account  (\ref{15.05.1}), we obtain
\begin{multline}\label{15.05.2}
a_1(y_1)(u_2+v_2)^2+b_1(y_1)(u_2+v_2)+c_1(y_1)-\sigma (u_2+y_1/2-v_2)^2\\ =
-\sigma(u_2-v_2)^2+a_2(y_1)(u_2-y_1/2+v_2)^2+b_2(y_1)(u_2-y_1/2+v_2)\\+c_2(y_1)+2\pi in(y_1),  \ \ |y_1|<\delta, \  u_2, v_2\in \mathbb{C},
\end{multline}
where the function $n(y_1)$ takes integer values. Considering the left-hand side and the right-hand side of  (\ref{15.05.2})  as polynomials in $u_2$ and $v_2$  and equating in  (\ref{15.05.2})  the coefficients of $u_2$ and $v_2$, we get
$$
b_1(y_1)-\sigma y_1=-a_2(y_1)y_1+b_2(y_1), \quad  b_1(y_1)+\sigma y_1=-a_2(y_1)y_1+b_2(y_1),\quad  |y_1|<\delta.
$$
This implies that $\sigma=0$. Hence  $\hat\nu_1(0, y_2)=\hat\nu_2(0, y_2)=1$ for all $y_2\in \mathbb{R}$.

3. We will prove the lemma by induction. Let $n=2$. We have, $$H=(I+\widetilde\alpha_2)(\mathbb{R}^2)=\{(0, y_2)\in \mathbb{R}^2: y_2\in \mathbb{R}\}.$$  Substituting $u_1=v_1=0$, $u_2=v_2=y_2/2$ into equation (\ref{14.05.8}) for the functions $\hat\mu_j(y_1, y_2)$,  we get that $\hat\mu_1(0, y_2)=\hat\mu_2(0, y_2)$ for all $y_2\in \mathbb{R}$.    It follows from what was proved in step 2 that $|\hat\mu_1(0, y_2)|=|\hat\mu_2(0, y_2)|=1$ for all $y_2\in \mathbb{R}$.   Hence  there is a real number   $x$ such that
$\hat\mu_j(0, y_2)=\exp\{ixy_2\}$ for all $y_2\in \mathbb{R}$, $j=1, 2$. Put $t_1=(x, -x)$, $t_2=(0, -x)$. Consider the distributions $\tau_j=\mu_j*E_{t_j}$. It is obvious that  $\hat\tau_j(0, y_2)=1$ for all $y_2\in \mathbb{R}$, $j=1, 2$. It follows from this that the distributions $\tau_j$ are supported in the annihilator $A(\mathbb{R}^2, H)$. It is obvious that $A(\mathbb{R}^2, H)=K$. Since $t_1+\alpha_2t_2=0$, the characteristic functions $\hat\tau_j(y)$   satisfy equation (\ref{11.04.1}).
By Lemma \ref{lem1}, this implies that if $\eta_j$ are independent random vectors with values in the space   $\mathbb{R}^2$ and distributions $\tau_j$, then the conditional distribution of the linear form   $M_2 = \eta_1 +\alpha_2\eta_2$ given $M_1 = \eta_1 + \eta_2$  is symmetric.  Since  $K={\Ker}~(I+\alpha_2)$, the restriction of the 
   operator    $\alpha_2$  to the subspace $K$ coincides with $-I$. It means that if we consider
$\eta_j$ as independent random vectors with values in    $K$, then the conditional distribution of the linear form   $M_2= \eta_1 -\eta_2$ given $M_1= \eta_1  + \eta_2$ is symmetric. It follows from Lemma \ref{lem3}, applying to the subspace $K$, that $\tau_1=\tau_2$.   Thus, when $n=2$ the lemma is proved.

Let $n>2$. We note that $\widetilde\alpha_n(y_1, y_2, \dots, y_n)=(-y_1, y_1-y_2, \dots, y_{n-1}-y_n)$.  Hence  equation (\ref{11.04.1}) for the functions $\hat\nu_j(y)$ takes the form
\begin{multline}\label{15.05.4}
\hat\nu_1(u_1+v_1, u_2+v_2,\dots ,u_n+v_n)\hat\nu_2(u_1-v_1, u_2+v_1-v_2,\dots, u_n+v_{n-1}-v_n)\\=
\hat\nu_1(u_1-v_1, u_2-v_2,\dots ,u_n-v_n)\hat\nu_2(u_1+v_1, u_2-v_1+v_2,\dots, u_n-v_{n-1}+v_n), \ \  u_j, v_j\in \mathbb{R}.
\end{multline}
Substituting  $u_1=u_2=\dots=u_{n-2}=0$, $v_1=v_2=\dots=v_{n-2}=0$  into (\ref{15.05.4})   we obtain
\begin{multline}\label{15.05.5}
\hat\nu_1(0,\dots, 0,  u_{n-1}+v_{n-1}, u_n+v_n)\hat\nu_2(0,\dots,0, u_{n-1}-v_{n-1}, u_n+v_{n-1}-v_n)\\=
\hat\nu_1(0,\dots, 0, u_{n-1}-v_{n-1}, u_n-v_n)\hat\nu_2(0,\dots,0, u_{n-1}+v_{n-1}, u_n-v_{n-1}+v_n), \ \  u_j, v_j\in \mathbb{R}.
\end{multline}
We see that   equation (\ref{15.05.5}), up to the notation, coincides with equation   (\ref{14.05.8}). Therefore, as proven in step 2, we have
$\hat\nu_1(0,\dots,0, y_n)=\hat\nu_2(0,\dots,0, y_n)=1$  for all $y_n\in \mathbb{R}$.
Hence  $|\hat\mu_1(0,\dots,0, y_n)|=|\hat\mu_2(0,\dots,0, y_n)|=1$ for all $y_n\in \mathbb{R}$. Put $L=\{(0,\dots,0, y_n)\in \mathbb{R}^n: y_n\in \mathbb{R}\}$.   Arguing as in the case  when   $n=2$ we can replace the distributions   $\mu_j$ by their shifts  $\tau_j$   in such a way that
the distributions $\tau_j$ are supported in the annihilator   $A(\mathbb{R}^n, L)=\mathbb{R}^{n-1}$. Moreover, if $\eta_j$ are independent random vectors with values in the space   $\mathbb{R}^n$ and distributions $\tau_j$, then the conditional distribution of the linear form   $M_2 = \eta_1 +\alpha_n\eta_2$ given 
$M_1 = \eta_1 + \eta_2$  is symmetric.
We note that $\alpha_n(\mathbb{R}^{n-1})=\mathbb{R}^{n-1}$  and the restriction of the  operator $\alpha_n$  to the subspace $\mathbb{R}^{n-1}$ coincides with the 
 operator $\alpha_{n-1}$.  It means that if we consider
$\eta_j$ as independent random vectors with values in     $\mathbb{R}^{n-1}$, then the conditional distribution of the linear form   $M_2= \eta_1 +\alpha_{n-1}\eta_2$ given $M_1= \eta_1+ \eta_2$ is symmetric.    The lemma is proved by induction.

The statement of the lemma can not be strengthened. Indeed, let $\omega$ be an arbitrary distribution supported in the subspace $K={\Ker}~(I+\alpha_n)$.
Let $t_1$ and $t_2$ be some vectors in $\mathbb{R}^n$ such that
$t_1+\alpha_n t_2=0$. Put $\mu_j=\omega*E_{t_j}$, $j=1, 2$.
Taking into account that the restriction of the  operator    $\alpha_n$  to the subspace $K$ coincides with $-I$,  Lemma \ref{lem1} and Lemma \ref{lem3} imply that if   $\xi_1$ and $\xi_2$ are independent   random vectors with values in  the space    $\mathbb{R}^n$ and distributions $\mu_j$, then
the conditional distribution of the linear form   $L_2 = \xi_1 + \alpha_n\xi_2$ given $L_1 = \xi_1 + \xi_2$ is symmetric. \hfill$\Box$  

\bigskip

The following statement implies from  Lemma   \ref{lem13}.
\begin{corollary}\label{co1}
Let    $\alpha$ be an invertible linear operator in the space $\mathbb{R}^n$.  Put $K={\Ker}~(I+\alpha)$ and suppose   $K\ne\{0\}$.  Let $\xi_1$ and $\xi_2$ be
independent random vectors with values in $\mathbb{R}^n$ and distributions
$\mu_1$ and $\mu_2$. Assume that   the conditional  distribution of the linear form  $L_2 = \xi_1 + \alpha\xi_2$ given $L_1 = \xi_1 + \xi_2$  is symmetric. Then there exists an  $\alpha$-invariant subspace $G$ satisfying the condition    $K\cap G=\{0\}$ and such that some shifts $\tau_j$  of the distributions $\mu_j$   are supported in the subspace $K\times G$. Moreover, if $\eta_j$ are independent   random vectors with values in the space $\mathbb{R}^n$ and distributions  $\tau_j$, then the conditional  distribution of the linear form
$M_2 = \eta_1 +\alpha\eta_2$ given $M_1 = \eta_1 + \eta_2$  is symmetric, and   the restriction of the  operator $\alpha$ to the subspace   $K\times G$ is of the form $(-I, \alpha_G)$.
\end{corollary}
{\bf Proof}. Represent the space $\mathbb{R}^n$ as a direct sum of two $\alpha$-invariant subspaces   $\mathbb{R}^n=F\times G$, where   $F$ is the root subspace corresponding  to the eigenvalue  $\lambda=-1$ of the  operator $\alpha$, and the  operator $I+\alpha$ is invertible in the subspace $G$. The  operator $\alpha$ can be written in the form  $\alpha=(\alpha_F, \alpha_G)$. In order not to complicate the notation  we assume that
$F=\mathbb{R}^p$, $G=\mathbb{R}^q$,
a basis in the space $\mathbb{R}^p$ is chosen in such a way that   $\mathbb{R}^p=\mathbb{R}^{n_1}\times\mathbb{R}^{n_2}\times\dots\times\mathbb{R}^{n_k}$, and the  operator      $\alpha_F$  is of the form  $\alpha_F=(\alpha_{n_1}, \alpha_{n_2},\dots, \alpha_{n_k})$, where   $\alpha_{n_j}=-I$, if $n_j=1$, and $\alpha_{n_j}$ is an invertible linear operator in the space $\mathbb{R}^{n_j}$  which is given by formula  $(\ref{16.05.1})$, if $n_j\ge 2$.   By Lemma \ref{lem1}, the characteristic functions  $\hat\mu_j(y)$   satisfy equation  (\ref{11.04.1}) which takes the form
\begin{multline}\label{16.05.2}
\hat\mu_1(u_1+v_1, u_2+v_2)\hat\mu_2(u_1+\widetilde\alpha_F v_1, u_2+\widetilde\alpha_G v_2)\\=
\hat\mu_1(u_1-v_1, u_2-v_2)\hat\mu_2(u_1-\widetilde\alpha_F v_1, u_2-\widetilde\alpha_G v_2), \quad  u_1, v_1\in \mathbb{R}^p, \ u_2, v_2\in \mathbb{R}^q.
\end{multline}
Substitute $u_2=v_2=0$ into (\ref{16.05.2}). Taking into account Lemma \ref{lem1}  and applying successively Lemma \ref{lem13} to each of the subspaces   $\mathbb{R}^{n_j}$, where    $n_j\ge 2$, we find as a result  from the obtained equation  that there exist vectors   $t_j\in \mathbb{R}^{p}$ such that the distributions     $\tau_j=\mu_j*E_{t_j}$ are supported in the subspace $K\times G$. Moreover, if $\eta_j$ are independent random vectors with values in $K\times G$  and distributions $\tau_j$, then the conditional distribution of the linear form   $M_2 = \eta_1 +\alpha\eta_2$ given $M_1 = \eta_1 + \eta_2$  is symmetric. Obviously, the restriction of the  operator $\alpha$ to the subspace   $K\times G$ is of the form $(-I, \alpha_G)$. \hfill$\Box$
\begin{remark}\label{re1}
Assume that under the conditions of Theorem  \ref{th1}  $K\ne \{0\}$, i.e. $\lambda=-1$ is an eigenvalue of the operator  $\alpha$.  Assume also that the root subspace corresponding to the eigenvalue $\lambda=-1$ does not coincide with the eigenspace. Then, as was proved in Corollary \ref{co1}, some shifts of the distributions
   $\mu_j$  are supported in a proper subspace of the space $\mathbb{R}^{n}$. Hence     $\mu_j$ are singular distributions.
\end{remark}
{\bf Proof of Theorem \ref{th1}.}  By Theorem B, if  $K=\{0\}$, i.e. $\lambda=-1$ is not an eigenvalue of the  operator $\alpha$, then $\mu_j$ are Gaussian distributions. Therefore, we   assume   $K\ne\{0\}$, i.e.   $\lambda=-1$ is an eigenvalue of the  operator $\alpha$. Let $F$ be the root subspace corresponding  to the eigenvalue  $\lambda=-1$ of the  operator $\alpha$. Corollary \ref{co1}
allows us to prove the theorem, assuming that $\alpha_F=-I$.   In other words, the root subspace corresponding to the eigenvalue $\lambda=-1$ of the  operator $\alpha$ is the eigenspace.  Then equation (\ref{16.05.2}) takes the form
\begin{multline}\label{16.05.3}
\hat\mu_1(u_1+v_1, u_2+v_2)\hat\mu_2(u_1- v_1, u_2+\widetilde\alpha_G v_2)\\=
\hat\mu_1(u_1-v_1, u_2-v_2)\hat\mu_2(u_1+ v_1, u_2-\widetilde\alpha_G v_2), \quad u_1, v_1\in \mathbb{R}^p, \ u_2, v_2\in \mathbb{R}^q.
\end{multline}

Substitute $u_1=v_1=0$ into  (\ref{16.05.3}). Taking into account Lemma \ref{lem1},  it follows from   Theorem B     applying to the space   $\mathbb{R}^{q}$ that   $\hat\mu_j(0, y_2)$ are the characteristic functions of Gaussian distributions, i.e.
\begin{equation}\label{16.05.4}
\hat\mu_j(0, y_2)=\exp\{-\langle A_j y_2, y_2\rangle+i\langle b_j, y_2\rangle\}, \quad y_2\in \mathbb{R}^q,
\end{equation}
where $A_j$  is a symmetric positive semidefinite $q\times q$ matrix, $b_j\in \mathbb{R}^q$, $j=1, 2$. Substituting $u_2=v_2=0$ into   (\ref{16.05.3}) and taking into account (\ref{16.05.4}) we obtain that $b_1+\alpha_Gb_2=0$. Hence  we can replace the distributions $\mu_j$ by their shifts   $\tau_j=\mu_j*E_{-b_j}$   and suppose that $b_1=b_2=0$ in (\ref{16.05.4}), i.e.
\begin{equation}\label{16.05.5}
\hat\mu_j(0, y_2)=\exp\{-\langle A_j y_2, y_2\rangle\}, \quad y_2\in \mathbb{R}^q, \   j=1, 2.
\end{equation}
By Lemma \ref{lem5}, it follows from (\ref{16.05.5}) that for each fixed $y_1\in \mathbb{R}^p$ the function $\hat\mu_j(y_1, y_2)$ can be extended  to $\mathbb{C}^q$ as an entire function in  $y_2$. It is obvious that equation (\ref{16.05.3}) remains valid  for all $u_1, v_1  \in \mathbb{R}^p$, $u_2, v_2  \in \mathbb{C}^q$.

Substituting $u_1=v_1=y_1/2$, $u_2=v_2=0$ into (\ref{16.05.3}) we see that $\hat\mu_1(y_1, 0)=\hat\mu_2(y_1, 0)$ for all $y_1\in \mathbb{R}^p$. Make sure that if for a fixed $\widetilde y_1\in \mathbb{R}^p$ the inequalities
\begin{equation}\label{20.05.1}
\hat\mu_j(\widetilde y_1, 0)\ne 0, \quad  j=1, 2,
\end{equation}
hold, then $\hat\mu_j(\widetilde y_1, y_2)\ne 0$   for   all $y_2  \in \mathbb{C}^q$, $j=1, 2$. Suppose the contrary, let
$\hat\mu_1(\widetilde y_1, \widetilde y_2)=0$ for some $\widetilde y_2  \in \mathbb{C}^q$. Substitute $u_1=v_1=\widetilde y_1/2$, $u_2=(I+\widetilde\alpha_G)^{-1} \widetilde\alpha_G\widetilde y_2$, $v_2=(I+\widetilde\alpha_G)^{-1}\widetilde y_2$ into (\ref{16.05.3}). Then the left-hand side of equation   (\ref{16.05.3}) is equal to zero. In view of  (\ref{16.05.5}) and (\ref{20.05.1}),  the right-hand side of equation     (\ref{16.05.3}) is nonzero. If $\hat\mu_2(\widetilde y_1, \widetilde y_2)=0$ for some $\widetilde y_2  \in \mathbb{C}^q$, then substituting $u_1=-v_1= \widetilde  y_1/2$, $u_2=v_2=(I+\widetilde\alpha_G)^{-1}\widetilde y_2$ into (\ref{16.05.3}), we get the contradiction,  because   the left-hand side of equation   (\ref{16.05.3}) is equal to zero  and the right-hand side is not. Thus, if   inequalities (\ref{20.05.1}) hold,  then the function $\hat\mu_2(\widetilde y_1, y_2)$ is also nonzero for  all   $y_2  \in \mathbb{C}^q$.

So, we proved that if $\hat\mu_j(y_1, 0)\ne 0$, $j=1, 2$, for some $y_1\in \mathbb{R}^p$, then the functions $\hat\mu_j(y_1, y_2)$ can be extended  to $\mathbb{C}^q$   as entire functions in $y_2$  without zeros. Hence  we have the representations
$$
\hat\mu_j(y_1, y_2)=\exp\{Q_j(y_1, y_2)\}, \quad  y_2  \in \mathbb{C}^q, \ j=1, 2,
$$
where  $Q_j(y_1, y_2)$ are  entire functions in   $y_2$ in $\mathbb{C}^q$.

By Lemma \ref{lem5}, it follows from inequality   (\ref{16.05.6}) and (\ref{16.05.5}) that  by   Hadamard's theorem on the representation of entire functions of finite order,    the restriction of the functions  $Q_j(y_1, y_2)$
to each complex plane in  $\mathbb{C}^q$ passing through zero  are polynomials of degree at most 2. Hence  the functions   $Q_j(y_1, y_2)$ are   polynomials 
of degree at most 2 in $y_2$. Thus, we have a representation
\begin{equation}\label{16.05.7}
\hat\mu_j(y_1, y_2)=\exp\{\langle A_j(y_1) y_2, y_2\rangle+\langle b_j(y_1), y_2\rangle+c_j(y_1)\}, \quad y_2\in \mathbb{C}^q,
\end{equation}
where $A_j(y_1)$ are symmetric complex $ q\times q $ matrices, $b_j(y_1)\in \mathbb{C}^q$, $c_j(y_1)\in \mathbb{C}$.

Assume   $\hat\mu_j(y_1, 0)\ne 0$, $j=1, 2$. Substituting  $u_1=v_1=y_1/2$ into  (\ref{16.05.3})  and taking into account (\ref{16.05.7}), we find from the obtained equation
\begin{multline}\label{16.05.8}
\langle A_1(y_1)(u_2+v_2), u_2+v_2\rangle+\langle A_2(0)(u_2+\widetilde\alpha_Gv_2), u_2+\widetilde\alpha_Gv_2\rangle\\=\langle A_1(0)(u_2-v_2), u_2-v_2\rangle+\langle A_2(y_1)(u_2-\widetilde\alpha_Gv_2), u_2-\widetilde\alpha_Gv_2\rangle,\quad u_2, v_2\in \mathbb{C}^q.
\end{multline}
\begin{equation}\label{17.05.3}
\langle b_1(y_1), u_2+v_2\rangle+\langle b_2(0), u_2+\widetilde\alpha_Gv_2\rangle=\langle b_1(0), u_2-v_2\rangle+\langle b_2(y_1), u_2-\widetilde\alpha_Gv_2\rangle,\quad u_2, v_2\in \mathbb{C}^q.
\end{equation}
The equalities
\begin{equation}\label{17.05.2}
A_1(y_1)+A_2(0)=A_1(0)+A_2(y_1)
\end{equation}
and
\begin{equation}\label{16.05.9}
A_1(y_1)+\alpha_GA_2(0)\widetilde\alpha_G=
A_1(0)+\alpha_GA_2(y_1)\widetilde\alpha_G
\end{equation}
follow from  (\ref{16.05.8}). Substituting  $u_1=y_1$,  $v_1=0$ into  (\ref{16.05.3}) and taking into account (\ref{16.05.7}), we get from the obtained equation
\begin{multline}\label{16.05.10}
\langle A_1(y_1)(u_2+v_2), u_2+v_2\rangle+\langle A_2(y_1)(u_2+\widetilde\alpha_Gv_2), u_2+\widetilde\alpha_Gv_2\rangle\\=\langle A_1(y_1)(u_2-v_2), u_2-v_2\rangle+\langle A_2(y_1)(u_2-\widetilde\alpha_Gv_2), u_2-\widetilde\alpha_Gv_2\rangle,\quad u_2, v_2\in \mathbb{C}^q.
\end{multline}
Equation (\ref{16.05.10}) implies that
\begin{equation}\label{16.05.11}
A_1(y_1)+\alpha_GA_2(y_1)=0.
\end{equation}
Taking into account that  $I+\alpha_G$ is an invertible   operator,   (\ref{16.05.5}), (\ref{16.05.9}) and (\ref{16.05.11}) imply that $A_1(y_1)=A_1(0)=-A_1$. Then it follows from  (\ref{16.05.5}) and (\ref{17.05.2}) that $A_2(y_1)=A_2(0)=-A_2$.
Thus, we have proved that if $\hat\mu_j(y_1, 0)\ne 0$, then
\begin{equation}\label{17.05.5}
A_1(y_1)=-A_1, \quad A_2(y_1)=-A_2.
\end{equation}
We find from (\ref{17.05.3}) that
\begin{equation}\label{17.05.4}
b_1(y_1)+b_2(0)=b_1(0)+b_2(y_1), \quad b_1(y_1)+\alpha_Gb_2(0)=-b_1(0)-\alpha_Gb_2(y_1).
\end{equation}
It follows from (\ref{16.05.5}) and (\ref{16.05.7}) that $b_1(0)=b_2(0)=0$. Given this, and taking into account that $I+\alpha_G$ is an invertible 
operator,   (\ref{17.05.4})  implies that if $\hat\mu_j(y_1, 0)\ne 0$, then
$b_1(y_1)=b_2(y_1)=0$ 
and hence  in view of (\ref{16.05.7}) and  (\ref{17.05.5}), the  representations
\begin{equation}\label{20.05.3}
\hat\mu_j(y_1, y_2)=\exp\{-\langle A_j y_2, y_2\rangle+c_j(y_1)\}, \quad y_2\in \mathbb{C}^q, \quad j=1, 2,
\end{equation}
are valid for the functions $\hat\mu_j(y_1, y_2)$. Moreover,   in this case substituting  $u_1=v_1=y_1/2$,  $u_2=v_2=0$ in equation  (\ref{16.05.3}) and taking into account  (\ref{20.05.3}),   we obtain that
\begin{equation}\label{17.05.7}
\hat\mu_1(y_1,0)=\hat\mu_2(y_1,0)=\exp\{c_1(y_1)\}=\exp\{c_2(y_1)\},\ \ y_1\in \mathbb{R}^p.
\end{equation}

Assume   $\hat\mu_j(\widetilde y_1, 0)=0$, $j=1, 2$, for some   $\widetilde y_1\in \mathbb{R}^p$. We verify that then $\hat\mu_j(\widetilde y_1, y_2)=0$    for  all $y_2  \in \mathbb{C}^q$, $j=1, 2$.   Put $u_1=v_1=\widetilde y_1/2$,  $u_2=-v_2=y_2$ in equation (\ref{16.05.3}). Then the left-hand side of equation    (\ref{16.05.3}) is equal to zero. In view of (\ref{16.05.5}), we have $\hat\mu_2(\widetilde y_1, (I+\widetilde\alpha_G)y_2)=0$ for all  $y_2  \in \mathbb{C}^q$. 
Since $(I+\widetilde\alpha_G)$ is an invertible  operator, we have $\hat\mu_2(\widetilde y_1, y_2)=0$ for all  $y_2  \in \mathbb{C}^q$. Substituting $u_1=-v_1=\widetilde y_1/2$,  $u_2=\widetilde\alpha_Gy_2$, $v_2=-y_2$ into equation (\ref{16.05.3}), we make sure that  $\hat\mu_1(\widetilde y_1, y_2)=0$ for all   $y_2  \in \mathbb{C}^q$.

Denote by $\gamma_j$ the symmetric Gaussian distribution in the space $\mathbb{R}^q$ with the characteristic function
\begin{equation}\label{19.05.2}
\hat\gamma_j(y_2)=\exp\{-\langle A_jy_2, y_2\rangle\}, \quad y_2\in \mathbb{R}^q, \ j=1, 2.
\end{equation}
In view of  $\hat\mu_1(y_1, 0)=\hat\mu_2(y_1, 0)$ for all $y_1\in \mathbb{R}^p$, denote by   $\omega$   a distribution in the space   $\mathbb{R}^p$ with the characteristic function
\begin{equation}\label{20.05.5}
\hat\omega(y_1)=\hat\mu_1(y_1,0)=\hat\mu_2(y_1,0), \quad y_1\in\mathbb{R}^p.
\end{equation}
Taking into account that $\hat\omega(y_1)=0$ if and only if   $\hat\mu_2(y_1, y_2)=0$ for  all $y_2  \in \mathbb{R}^q$, (\ref{20.05.3})--(\ref{20.05.5}) imply the representations
$$
\hat\mu_j(y_1, y_2)=\exp\{-\langle A_jy_2, y_2\rangle\}\hat\omega(y_1), \quad y_1\in \mathbb{R}^p, \ y_2\in \mathbb{R}^q, \ j=1, 2.
$$
This implies the statement of the theorem.

Note that we also proved that the intersection of each of the 
supports of the symmetric Gaussian distributions    
$\gamma_j$ with the root subspace corresponding  to the eigenvalue  
$\lambda=-1$ of the  operator $\alpha$ is equal to zero.

The statement of the theorem can not be strengthened. Indeed, consider     an invertible linear operator $\alpha$ in the space $\mathbb{R}^n$. Let $G$ be an
$\alpha$-invariant subspace in $\mathbb{R}^n$ such that the 
 operator  $I+\alpha$ is invertible in $G$. Let $G$ be isomorphic to   $\mathbb{R}^q$. Denote by $\gamma_j$ a Gaussian distribution in the space $\mathbb{R}^q$ with the characteristic function  (\ref{19.05.2}). Moreover, we will also assume that  the equality
\begin{equation}\label{19.05.3}
A_1+A_2\widetilde\alpha_G=0
\end{equation}
holds. Put $K={\Ker}~(I+\alpha)$. Let $\omega$ be a distribution on   $K$.
Let $x_1$ and $x_2$ be some vectors in $\mathbb{R}^n$ such that
\begin{equation}\label{13.07.1}
x_1+\alpha x_2=0.
\end{equation}
Put $\mu_j=\gamma_j*\omega*E_{x_j}$, $j=1, 2$. Let $\xi_1$ and $\xi_2$ be
independent random vectors with values in $\mathbb{R}^n$ and distributions
$\mu_j$.
It is easy to  see that (\ref{19.05.3}) implies that the characteristic functions  $\hat\gamma_j(y)$ satisfy equation (\ref{11.04.1}). It follows from Lemma \ref{lem1} and Lemma \ref{lem3}, applying to the subspace  $K$, that  the characteristic function $\hat\omega(y)$ satisfies equation (\ref{11.04.1}).  Moreover, (\ref{13.07.1}) implies that the functions  $(x_j, y)$  satisfy equation (\ref{11.04.1}).  Hence  the characteristic functions  $\hat\mu_j(y)$ satisfy equation (\ref{11.04.1}).  By Lemma \ref{lem1},  the conditional  distribution of the linear form  $L_2 = \xi_1 + \alpha\xi_2$ given $L_1 = \xi_1 + \xi_2$  is symmetric.
\hfill$\Box$

\bigskip

We complement Theorem \ref{th1} by  a more detailed description of   possible distributions $\mu_j$ in Theorem \ref{th1} for the   space   $\mathbb{R}^2$ depending on the spectrum of the linear operator $\alpha$. If $\lambda$ is an eigenvalue of the  operator  $\alpha$, denote by $L_{\lambda}$ the corresponding eigenspace. In particular, if $K={\Ker}~(I+\alpha)\ne \{0\}$, then $K=L_{-1}$. We consider two cases: $\lambda=-1$ is an eigenvalue of the operator $\alpha$ or not.

$1$.  $\lambda=-1$ is an eigenvalue of the operator $\alpha$. Then the characteristic equation of the operator  $\alpha$ has only real   roots. Two cases are possible.

1A. The  operator $\alpha$ has another eigenvalue $\lambda=\lambda_0$, where $\lambda_0\ne -1$. In this case in a basis consisting of eigenvectors of $\alpha$, the diagonal matrix $\alpha={\rm diag}\{-1, \lambda_0\}$ corresponds to the operator $\alpha$.   Theorem \ref{th1} easily implies the following alternative. If $\lambda_0> 0$, then $\mu_j=\omega*E_{x_j}$, where $\omega$ is a distribution supported in    $K$, $j=1, 2$. In so doing  $K$ is a one-dimensional subspace of   $\mathbb{R}^2$. If $\lambda_0< 0$, then $\mu_j=\gamma_j*\omega$, where $\gamma_j$ are Gaussian distributions supported in    $L_{\lambda_0}$, and $\omega$ is a distribution supported in $K$.

1B. The operator $\alpha$ has the only eigenvalue   $\lambda=-1$. If the root subspace corresponding  to the eigenvalue   $\lambda=-1$  does not coincide with $K$, then the Jordan cell
$$
\alpha=\left(
\begin{array}{cc}
-1 & 1\\
0 & -1
\end{array}
\right)
$$
in a suitable basis corresponds to the operator $\alpha$.
In this case  by Lemma \ref{lem13},
$\mu_j=\omega*E_{x_j}$, where $\omega$ is a distribution supported in  $K$, $j=1, 2$. In so doing  $K$ is a one-dimensional subspace of   $\mathbb{R}^2$. If the root subspace corresponding  to the eigenvalue   $\lambda=-1$   coincides with $K$, then    $\alpha=-I$, and by Lemma \ref{lem3}, $\mu_1=\mu_2=\mu$, where $\mu$ is an arbitrary distribution on $\mathbb{R}^2$.

$2$.   $\lambda=-1$ is not  an eigenvalue of the operator $\alpha$.   By Theorem B, in this case $\mu_j$ are Gaussian distributions.  Find out what can be said about the supports of the Gaussian distributions   $\mu_j$, in particular,  when $\mu_j$ are degenerate distributions.

The characteristic functions $\hat\mu_j(y_1, y_2)$ are represented in the form
\begin{equation}\label{27.05.1}
\hat\mu_j(y_1, y_2)=\exp\{-\langle A_j (y_1, y_2),  (y_1, y_2)\rangle+i\langle b_j,(y_1, y_2)\rangle\}, \quad (y_1, y_2)\in \mathbb{R}^2,
\end{equation}
where $A_j$  is a symmetric positive semidefinite $2\times 2$ matrix, $b_j\in \mathbb{R}^2$, $j=1, 2$. By Lemma \ref{lem1}, the characteristic functions  $\hat\mu_j(y_1, y_2)$ satisfy equation  (\ref{11.04.1}). Depending on the spectrum of the operator    $\alpha$,  a matrix of a rather simple form  in a suitable basis corresponds to the operator $\alpha$. Substituting (\ref{27.05.1}) into (\ref{11.04.1}) we get that the matrices $A_j$ satisfy the equation
\begin{equation}\label{27.05.2}
A_1+A_2\widetilde\alpha=0.
\end{equation}
The description of solutions of equation  (\ref{27.05.2}) in the class of
symmetric positive semidefinite $2\times 2$ matrices implies the description of the supports of the  Gaussian distributions $\mu_j$.
 Two cases are possible.

2A. The characteristic equation of the operator  $\alpha$ has only real   roots.

(i). The operator $\alpha$  has two different eigenvalues    $\lambda=\lambda_1$ and $\lambda=\lambda_2$.
In this case in a basis consisting of eigenvectors of $\alpha$, the diagonal matrix $\alpha={\rm diag}\{\lambda_1, \lambda_2\}$ corresponds to the 
operator $\alpha$. Then the following statements are valid.  If $\lambda_1>0$ and $\lambda_2>0$, then $\mu_j$  are degenerate distributions. If $\lambda_1\lambda_2<0$, then $\mu_j=\gamma_j*E_{x_j}$, where $\gamma_j$ are symmetric Gaussian distributions supported in   $L_{\lambda_j}$, where $\lambda_j<0$. If $\lambda_1<0$ and $\lambda_2<0$, then $\mu_j=\gamma_j*E_{x_j}$, where $\gamma_j$ are symmetric Gaussian distributions. In so doing  either  $\gamma_j$ are degenerate distributions concentrated at zero,  or the  supports of $\gamma_j$ are   one of subspaces   $L_{\lambda_j}$, or   the supports of  $\gamma_j$ are the space $\mathbb{R}^2$.

(ii). The operator $\alpha$ has the only eigenvalue  $\lambda=\lambda_0$. If the root space corresponding to the eigenvalue $\lambda=\lambda_0$ does not coincide  with $L_{\lambda_0}$,
then the Jordan cell
$$\alpha=\left(
\begin{array}{cc}
\lambda_0 & 1\\
0 & \lambda_0
\end{array}
\right)
$$
in a suitable basis corresponds to the  operator $\alpha$. Then, if  $\lambda_0>0$, then $\mu_j$ are degenerate distributions. If $\lambda_0<0$, then   $\mu_j=\gamma_j*E_{x_j}$, where $\gamma_j$ are symmetric Gaussian distributions supported in   $L_{\lambda_0}$.
If the root space corresponding to the eigenvalue $\lambda=\lambda_0$  coincides with $L_{\lambda_0}$, then $\alpha=\lambda_0 I$.  In this case, if  $\lambda_0>0$, then $\mu_j$ are degenerate distributions. If $\lambda_0<0$, then   $\mu_j=\gamma_j*E_{x_j}$, where $\gamma_j$ are symmetric Gaussian distributions. In so doing  either  $\gamma_j$ are degenerate distributions concentrated at zero,  or the  supports of $\gamma_j$ are the same  one-dimensional subspace  of  $\mathbb{R}^2$, or  the supports of $\gamma_j$ are the space $\mathbb{R}^2$.

2B.  The roots of the characteristic equation of the operator  $\alpha$ are of the form $\lambda_1=x_0+iy_0$, $\lambda_1=x_0-iy_0$, where $y_0\ne 0$. Then  the matrix
$$\alpha=\left(
\begin{array}{cc}
x_0 & y_0\\
-y_0 & x_0
\end{array}
\right)
$$
in a suitable basis corresponds to the  operator $\alpha$. In this case $\mu_j$  are degenerate distributions.

 \medskip
 
{\bf Acknowledgements}

\medskip

I thank the reviewer for getting my attention to articles \cite{Ai}--\cite{Al2}.

\medskip

This article was written during my stay at 
the Department of Mathematics University of Toronto as a Visiting Professor. 
I am very grateful to Ilia Binder for his invitation and support.


\end{document}